\documentclass[11pt,a4paper]{article}
\usepackage{amsmath,amssymb}
\newcommand{\e}{\varepsilon}

\newcommand{\D}{\Delta}

\newcommand{\n}{\nabla}
\newcommand{\N}{\frac{N}{2}}
\newcommand{\NN}{\frac{N}{p}}

\newcommand{\p}{\partial}

\newcommand{\R}{\mathbb{R}}

\newcommand{\de}{\delta}

\newtheorem{definition}{Definition}
\newtheorem{theorem}{Theorem}
\newtheorem{proposition}{Proposition}
\newtheorem{corollaire}{Corollary}

\newtheorem{remarka}{Remark}

\newtheorem{lemme}{Lemma}

\title{Existence of global strong solutions for the shallow-water equations with large initial data}
\author{Boris Haspot \thanks{Ceremade UMR CNRS 7534
Universit\'e de Paris IX- Dauphine,
Place du MarŽchal DeLattre De Tassigny
75775 PARIS CEDEX 16 , haspot@ceremade.dauphine.fr }}

\date{}
\begin{document}

\maketitle
\begin{abstract}
This work is devoted to the study of a viscous shallow-water system with friction and capillarity term. We prove in this paper the existence of global strong solutions for this system with some choice of large initial data when $N\geq 2$ in critical spaces for the scaling of the equations. More precisely, we introduce as in \cite{Hprepa} a new unknown,\textit{a effective velocity} $v=u+\mu\n\ln h$ ($u$ is the classical velocity and $h$ the depth variation of the fluid) with $\mu$ the viscosity coefficient which simplifies the system and allow us to cancel out the coupling between the velocity $u$ and the depth variation $h$. We obtain then the existence of global strong solution if $m_{0}=h_{0}v_{0}$ is small in $B^{\N-1}_{2,1}$ and $(h_{0}-1)$  large in $B^{\N}_{2,1}$. In particular it implies that the classical momentum $m_{0}^{'}=h_{0} u_{0}$ can be large in $B^{\N-1}_{2,1}$, but small when we project $m_{0}^{'}$ on the divergence field. These solutions are in some sense \textit{purely compressible}. We would like to point out that the friction term  term has a fundamental role in our work inasmuch as coupling with the pressure term it creates a damping effect on the effective velocity.
\end{abstract}
\section{Introduction}
We consider the viscous shallow water model with friction and capillarity term. This model is also called by the french community the Saint-Venant equations and is generally used in oceanography. Indeed it allows to model vertically averaged flows  in terms of the horizontal mean velocity field $u$ and the depth variation $h$ due to the free surface. In the rotating framework, the model is described by the following system:
\begin{equation}
\begin{cases}
\begin{aligned}
&\frac{\p}{\p t}h+{\rm div}(h u)=0,\\
&\frac{\p}{\p t}(h u)+{\rm div}(h
u\otimes u)-\rm div(\mu h\, Du)+\frac{\n h}{Fr^{2}}+r h u={\rm div}K,
\end{aligned}
\end{cases}
\label{3systeme}
\end{equation}
where ${\rm div}K$ is the free surface tension tensor which reads as follows:
\begin{equation}
{\rm div}K
=\n\big(h\kappa(h)\D h+\frac{1}{2}(\kappa(h)+h\kappa^{'}(h))|\n h|^{2}\big)
-{\rm div}\big(\kappa(h)\n h\otimes\n h\big).
\label{divK}
\end{equation}
$\kappa$ is the coefficient of free surface tension and is a regular function of the form $\kappa(h)=\kappa h^{\alpha}$ with $\alpha\in\R$. In the sequel we will assume that $\alpha=-1$.
$Fr > 0$ denotes the Froude number. System (\ref{3systeme}) is supplemented with initial conditions
\begin{equation}
h_{/t=0}= h_{0}, (hu)_{/t=0}=m^{'}_{0}.
\label{data}
\end{equation}
This model is derived from the three-dimensional Navier-Stokes equations with free surface, where the normal stress is determined from the air pressure and capillary effects. The turbulent regime ($r\geq 0$) is obtained from the friction condition on the bottom, see \cite{P}. $\mu$ is the viscosity coefficient and verifies $\mu>0$ and $D u=(\n u+^{t}\n u)/2$ is the strain tensor.\\
Several physical models arise as a particular case of system (\ref{3systeme}):
\begin{itemize}
\item when $\kappa=r=0$, (\ref{3systeme}) represents compressible Navier-Stokes model with shallow-water viscosity coefficients.
\item when $\kappa>0$ and $r=0$, then (\ref{3systeme}) describes the Korteweg system which models mixture liquid-vapor.
\end{itemize}
We would  like to point out also the theoretical aspect of system (\ref{3systeme}). Indeed in the case $r=0$ the system (\ref{3systeme}) corresponds  to the classical Korteweg system which models a liquid-vapour mixture. Let us mention that the Korteweg system also is used in a purely theoretical interest consisting in the selection of the physically relevant solutions of the Euler model by a vanishing capillarity-viscosity limit  (in particular when the system is not strictly hyperbolic which is typically the case when the pressure is Van der Waals).  Indeed in this case at least when $N=1$ it is not possible to apply the classical theory of Lax for the Riemann problem (see  \cite{Lax}) and of Glimm (see \cite{Glimm}) with small initial data in $BV$ in order to obtain the existence of global entropic solution (we refer also to the work of Bianchini and Bressan see \cite{BB} for the uniqueness).  In particular in this direction, recently in \cite{CH} with F. Charve we prove that the global strong solution of the Korteweg system in one dimension (we obtain also in this paper the existence of global strong solution in one dimension for Korteweg system)  converges in the setting of a $\gamma$ law for the pressure ($P(\rho)=a\rho^{\gamma}$, $\gamma>1$) to entropic solution of the compressible Euler equations. In particular it justifies that the Korteweg system is suitable for selecting the physical solutions in the case where the Euler system is strictly hyperbolic. The problem remains however open for a Van der Waals pressure.
\\
Now before investigating the problem of global strong solution for the system (\ref{3systeme}), we would like to recall the energy inequalities associated to this system and in particular describing some results about the  existence of global weak solutions for the system (\ref{3systeme}). Let $\bar{h}>0$ be a constant depth variation (in the sequel we will assume that $\bar{h}=1$),   and let $\Pi $ be defined
by:
$$\Pi(s)=s\biggl(\int^{s}_{\bar{h}}\frac{P(z)}{z^{2}}dz-\frac{P(\bar{h})}{\bar{h}}\biggl),$$
so that $P(s)=s\Pi^{'}(s)-\Pi(s)\, ,\,\Pi^{'}(\bar{h})=0$.
Multiplying the equation of momentum conservation in the system
(\ref{3systeme}) by $u$ and integrating by parts over $\R^{N}$,
we obtain the following
estimate:
\begin{equation}
\begin{aligned}
&\int_{\R^{N}}\big(\frac{1}{2}h
|u|^{2}+(\Pi(h)-\Pi(\bar{h}))+\frac{\kappa(h)}{2}|\nabla h|^{2}\big)(t)dx
+\frac{1}{2}\int_{0}^{t}\int_{\R^{N}}\mu h
|D(u)|^{2}dxdt\\
&\hspace{4cm}\leq\int_{\R^{N}}\big(\frac{|m^{'}_{0}|^{2}}{2 h}+(\Pi(h_{0})-\Pi(\bar{h}))
+\frac{\kappa(h_{0})}{2}|\nabla  h_{0}|^{2}\big)dx,
\label{3inegaliteenergie1}
\end{aligned}
\end{equation}
Here we can observe that if we assume that the initial data are such that:
\begin{equation}
\begin{aligned}
&h_{0}\ln(\frac{h_{0}}{e})\in L^{1},\;\;\sqrt{\kappa(h_{0})}\n h_{0}\in L^{2},\\
&\sqrt{h_{0}}|u_{0}|\in L^{2}(\R^{N}),
\end{aligned}
\label{initial}
\end{equation}
then we have the following estimates:
\begin{equation}
\begin{aligned}
&h|u|^{2}\in L^{\infty}(L^{1}),\;(\Pi(h)-\Pi(\bar{h}))\in L^{\infty}(L^{1}),\\
&\frac{\kappa(h)}{2}|\nabla h|^{2}\in L^{\infty}(L^{1}),\;Du\in L^{2}(L{2}).
\end{aligned}
\label{energie}
\end{equation}
One of the main difficulty in order to obtain the existence of global weak solution consists in dealing with the quadratic terms of the capillarity tensor. Indeed in order to have the stability of a sequence $(h_{n},u_{n})_{n\in\mathbb{N}}$ of global weak solution for the system (\ref{3systeme}), it is crucial to give a sense to the quadratic terms in the gradient of the depth variation $\kappa (h_{n})\n h_{n}\otimes\n h_{n}$ which are only uniformly bounded in $L^{\infty}(L^{1}(\R^{N}))$. In particular it implies only a convergence in the sense of the measure of  $\kappa (h_{n})\n h_{n}\otimes\n h_{n}$ which is not sufficient to conclude by standard compactness argument. That is why the problem of the existence of global weak solution remains open in the general case. 
Before stating our main result on the existence of global strong solution, we would like to recall what is known on the existence of global weak and strong solution in the different standard configurations. And in particular we would like to emphasize on the results obtained in \cite{Hprepa,Hprepa1} where we have discovered new entropies on the Korteweg system ($r=0$) for a specific choice on the capillarity $(\kappa(\rho)=\frac{\kappa}{\rho}$ with $\kappa>0$), allowing in particular to prove the existence of global weak solution. It is precisely with this type of capillarity that we will work with in the sequel.
\subsection*{Existence of global weak solutions}
\subsubsection*{Case $\kappa=0$ and $r=0$, the compressible Navier-Stokes system}
When the viscosity coefficients are constant and the pressure is a $\gamma$ law $P(h)=a h^{\gamma}$, with $a>0$ and
$\gamma>1$, P-L. Lions in \cite{1L2} proved the global existence of
variational solutions $(h,u)$ to (\ref{3systeme}) with $\kappa=r=0$
for $\gamma> \frac{N}{2}$ if $N\geq 4$, $\gamma\geq \frac{3N}{N+2}$
if $N=2,3$ and initial data $(\rho_{0},m_{0})$ such that:
$$\Pi(h_{0})-\Pi(\bar{h}),\;\;\frac{|m_{0}|^{2}}{h_{0}}\in
L^{1}(\R^{N}).$$ These solutions are weak solutions in the classical
sense for the equation of mass conservation and for
the equation of the momentum. Notice that the main difficulty for proving Lions' theorem consists
in exhibiting strong compactness properties of the height $h$ in
$L^{p}_{loc}$ spaces required to pass to the limit in the pressure
term $P(h)=a h^{\gamma}$.\\
Let us mention that Feireisl in \cite{fF}  generalized the result to
$\gamma>\frac{N}{2}$ in establishing that we can obtain renormalized
solution without imposing that $h\in L^{2}_{loc}$, for this he
introduces the concept of oscillation defect measure evaluating the
lost of compactness.\\
Concerning the shallow-water system when $\mu(h)=\mu h$, the main difficultywhen dealing with vanishing viscosity coefficients on vacuum is that the velocity cannot even be defined when the density vanishes. In particular we lose the information on $\n u$ in $L^{2}((0,T)\times\R^{N})$.
The main difficulty, to prove the stability of the solutions, is to pass to the limit in the term $\rho u\otimes u$ (which requires the strong convergence of $\sqrt{\rho}u$). Mellet and Vasseur  in \cite{fMV1} obtain the stability of global weak solution by using new entropies on the velocity and the density.
\subsubsection*{Case $r=0$, Korteweg system}
In the capillary case in contrast to the non capillary case, we can easily deal with the
pressure term in order to obtain stability results. However let us emphasize at this point that the energy estimates do not provide any $L^{\infty}$ control on the
density from below or from above. Indeed, even in dimension $N=2$,
$H^{1}$ functions are not necessarily locally bounded. Thus, vacuum
patches are likely to form in the fluid in spite of the presence of capillary forces, which are expected to smooth out the density. It explains why it is so difficult to obtaining the existence of global strong solution in dimension $N=2$ and even global weak solution. 
In \cite{fH2},we obtain the existence of global weak solution for specific choices of the capillary coefficients and with general viscosity coefficient but with small initial data in the energy space. More recently in \cite{Hprepa} and \cite{Hprepa1}, we prove the existence of global weak solution with large initial data for the Korteweg system (when $\kappa(h)=\frac{1}{h}$) with Saint-Venant viscosity coefficients. Indeed with this choice of capillarity coefficient we are able to obtain new entropy inequalities. To do this, we introduce a new unknown $v=u+\mu\n\ln\rho$ called \textit{effective velocity}, and we are able to show some  gain of integrability on $v$ which allow us to deal with the terms of the type $\rho u\otimes v$ (we refer to \cite{Hprepa1} for more details).
\subsubsection*{Existence of global strong solutions}
\subsubsection*{Case $\kappa=0$ and $r=0$, the compressible Navier-Stokes system}
We  refer to \cite{arma} for the existence of global strong solution with small initial data in critical space for the scaling of the equations. More precisely $(h_{0}-1,u_{0})$ belongs to $B^{\NN}_{p,1}\times B^{\frac{N}{p_{1}}-1}_{p_{1},1}$ for suitable chosen $p$ and $p_{1}$.
\subsubsection*{Case $r=0$, Korteweg system}
Let us mention briefly that the existence of global strong solutions in critical spaces for the scaling of the equations for $N\geq2$ is known since the works by  R. Danchin and B. Desjardins (when the viscosity coefficient and the capillary coefficients are constant, see also \cite{fH1}) in \cite{fDD} where the initial data $(h_{0}-1,h_{0}u_{0})$
belong to the Besov spaces $B^{\N}_{2,1}\times B^{\N-1}_{2,1}$ and are chosen small enough. In \cite{Hprepa}, we improve this result by working in a larger space of initial data, more precisely  $(h_{0}-1,h_{0}u_{0})$
is in the Besov space $B^{\N}_{2,\infty}\times B^{\N-1}_{2,\infty}$. In this paper for the same reason than in \cite{Hprepa1}, we are working with specific capillary coefficient $\kappa(h)=\frac{1}{h}$ because we are able to exhibit a specific structure on the unknown $\ln h$. In particular we are able to work with vortex initial data on the initial velocity in dimension $N=2$.
\subsection{Derivation of the model and non trivial explicit solutions}
The choice of system (\ref{3systeme}) is motivated by its energetic consistency, which has been stressed out from a physical point of view in \cite{13}. Compared with the Korteweg system ( see \cite{fK}) we also take into account the friction term $r h u$ physically justified to model the friction condition on the bottom of the ocean ( see \cite{15}). Naturally on a mathematical point of view,
the term $r h u$ does not add any difficulties for obtaining global weak solution or global strong solution with small initial data, however in our study this term shall turn out to be essential to ensure the existence of global strong solution with large initial data. Indeed coupled with the pressure term, he shall introduce a damping effect on a new unknown called \textit{effective velocity} and introduced in \cite{Hprepa,Hprepa1}. Roughly speaking this friction term allows to cancel out the coupling between the height $h$ and this effective velocity $v$, it is one of the main difficulty in order to obtain the existence of global strong solution for compressible Navier-Stokes equation (see \cite{CD, CMZ,arma}, in this case the coupling is between the velocity and the density). Indeed in this last case, it is difficult to obtain a damping effect on the density and the coupling between the velocity and the pressure terms impose a smallness condition on the initial density. In our case the fact that with the friction term we can rewrite the pressure term in some sense as a velocity, more precisely as \textit{the effective velocity} allows to avoid any condition of smallness on the initial depth variation $h_{0}$. We shall come back more in details on these considerations in the proof of theorem \ref{theo1}, we now would like to introduce this notion of \textit{effective velocity} which is one crucial tool for the proof of our result.\\
More precisely we are going to explain why with specific choice on the capillarity, Froude and friction   coefficient (as in  \cite{Hprepa} and \cite{Hprepa1}) we can exhibit a new structure onthe system via the introduction of this \textit{effective velocity} $v$. In particular in \cite{Hprepa,Hprepa1} we obtain new entropies which allows us toprove the existence of global weak solutions. 
We are now considering in the sequel the following physical coefficients:
\begin{equation}
\kappa(h)=\frac{\kappa}{h},\;\kappa=\mu^{2}\;\;\mbox{and}\;\;\frac{1}{Fr^{2}}=r\mu,
\label{hypothese}
\end{equation}
with $\mu>0$. By computation (see \cite{Hprepa}), we obtain the simplified system:
\begin{equation}
\begin{cases}
\begin{aligned}
&\p_{t}h+{\rm div}(h v)-\mu\D h=0,\\
&h\p_{t}v +h u\cdot\n v-\rm div(\mu h\,\n v)+r \,h\,v=0,
\end{aligned}
\end{cases}
\label{3systeme2}
\end{equation}
with $v=u+\mu\n\ln h$ the \textit{effective velocity}. For more details on the computation, we refer to \cite{Hprepa1}. When we write the system (\ref{3systeme2}) in function of the momentum $m=h v$, the system reads as follows:
\begin{equation}
\begin{cases}
\begin{aligned}
&\p_{t}h+{\rm div}m-\mu\D h=0,\\
&\p_{t}m +\rm div(\frac{m}{h}\otimes m)-\mu\D m+r\,m=0,
\end{aligned}
\end{cases}
\label{3system3}
\end{equation}
In particular we can observe that $m=0$ and $h_{1}$ such that:
\begin{equation}
\begin{cases}
\begin{aligned}
&\p_{t}h_{1}-\mu\D h_{1}=0,\\
&h_{1}(0,x)=h_{0}(x).
\end{aligned}
\end{cases}
\label{3system3}
\end{equation}
is a particular solution of (\ref{3system3}). When we consider the system (\ref{3systeme}), it means that $(h_{1},u_{1}=\mu\n\ln h_{1})$ is a non trivial solution. In particular this solution is \textit{purely compressible} as 
${\rm curl}u_{1}=0$.\\
In the sequel we are going to work around this non trivial solution $(h_{1},0)$ for the system (\ref{3system3}). If we now consider the system (\ref{3system3}), we can observe that the coupling between the momentum $m$ and the height $h$ via the pressure term has disappeared. As we explained previously the friction term and the pressure introduce here a damping effect via the new term $r\,m$. It is now possible to obtain the existence of  global strong solution when we assume only a condition of smallness on $m_{0}$ in order to deal with the non linear term ${\rm div}(\frac{m}{h}\otimes m)$ and in particular to use smallness argument for the term $\frac{m}{h}$.
\section*{Results}
Our main motivation concerns the existence of global strong solution with large initial data on the height $h_{0}$ but also on the initial velocity $u_{0}$. More precisely we will obtain global strong solution with a family of initial velocity $u_{0}$ such that the projection on gradient vector field  has large norm in $B^{\N-1}_{2,1}$.\\
We prove global well-posedness for system (\ref{3system3}
) in critical Besov space. To do this, we shall work around a constant depth variation $\bar{h}=1$, and to do this we introduce the following definition.
\begin{definition}
We set:
$$q_{0}=h_{0}-1.$$
\end{definition}
We can then rewrite the system (\ref{3system3}) in function of $(q,m)$, it gives:
\begin{equation}
\begin{cases}
\begin{aligned}
&\p_{t}q+{\rm div}m-\mu\D q=0,\\
&\p_{t}m +\rm div(\frac{m}{h}\otimes m)-\mu\D m+r\,m=0,
\end{aligned}
\end{cases}
\label{3system4}
\end{equation}
\begin{theorem}
Suppose that we are under the conditions (\ref{hypothese}). Assume that $m_{0}\in B^{\N-1}_{2,1}$ and $q_{0}\in B^{\N}_{2,1}$ with $h_{0}\geq c>0$.
Then there exists a constant $\e_{0}$ depending on $\frac{1}{h_{0}}$ such that if:
$$\|m_{0}\|_{B^{\N-1}_{2,1}}\leq\e_{0},$$
then there exists a unique global solution $(q,m)$ for system (\ref{3system4})
with $h$ bounded away from zero and,
$$
\begin{aligned}
&h\in \widetilde{C}(\R^{+},B^{\N}_{2,1}
)\cap L^{1}(\R^{+},
B^{\N+2}_{2,1})\;\;\;\mbox{and}\;\;\;\;m\in \widetilde{C}(\R^{+};B^{\N-1}_{2,1})
\cap L^{1}(\R{+},B^{\N-1}_{2,1}\cap B^{\N+1}_{2,1}).
\end{aligned}
$$
\label{theo1}
\end{theorem}
\begin{remarka}
In this theorem, for the first time up my knowledge we obtain a result on the existence of global strong solution for a compressible system without assuming any smallness hypothesis on the density as it is classically the case (see \cite{DG,CD,CMZ,arma}). Furthermore the physical moment $m^{'}_{0}=h_{0}u_{0}$ is large in $B^{\N-1}_{2,1}$, indeed we have:
$$m^{'}_{0}=m_{0}-\mu\n\rho_{0},$$
and here $m_{0}$ is small in $B^{\N-1}_{2,1}$ but $\n\rho_{0}$ could be arbitrary large. It means that for this family of large initial data ( when $m^{'}_{0}$ is the sum of a small momentum and of the gradient of the density), we have the existence of global strong solution. In particular the projection on the divergence field issmall in
$B^{\N-1}_{2,1}$ but the projection on the gradient field is large in $B^{\N-1}_{2,1}$. It is essentially related with the compressibility of the system.
\end{remarka}
\begin{remarka}
In the same way than in \cite{Hprepa}, we could consider the unknown $(\ln h,u)$ and in this case we could obtain the existence of global strong solution with small initial data in $(B^{\N}_{2,\infty},B^{\N-1}_{2,\infty})$ which is a larger space than in theorem \ref{theo1}. In particular it allows us to deal with the problem of vortex initial data, indeed in dimension $N=2$ we could choose initial velocity such that ${\rm curl}u_{0}$ is a bounded measure and in particular a Dirac $\delta_{0}$. However it would be not clear how to get global solution without smallness hypothesis on the initial density, indeed in the term $u\cdot\n v$ that we write under the form $v\cdot\n v-\mu\n\ln\rho\cdot\n v$ it is not clear how to deal with the quadratic term $\n\ln\rho\cdot\n v$. 
\end{remarka}
\begin{remarka}
Our method may be adapted to the $L^{p}$ framework (that is we now consider $q_{0}\in B^{\NN}_{p,1}$, $m_{0}\in  B^{\NN-1}_{p,1}$) when $1\leq p<+\infty$ as in \cite{fDD}.
\end{remarka}

Our paper is structured as follows. In section \ref{section2}, we give a few notation and briefly introduce the basic Fourier analysis
techniques needed to prove our result. In section \ref{section3}, we prove  the theorems\ref{theo1}.
\section{Littlewood-Paley theory and Besov spaces}
\label{section2}
Throughout the paper, $C$ stands for a constant whose exact meaning depends on the context. The notation $A\lesssim B$ means
that $A\leq CB$.
For all Banach space $X$, we denote by $C([0,T],X)$ the set of continuous functions on $[0,T]$ with values in $X$.
For $p\in[1,+\infty]$, the notation $L^{p}(0,T,X)$ or $L^{p}_{T}(X)$ stands for the set of measurable functions on $(0,T)$
with values in $X$ such that $t\rightarrow\|f(t)\|_{X}$ belongs to $L^{p}(0,T)$.
Littlewood-Paley decomposition  corresponds to a dyadic
decomposition  of the space in Fourier variables.
We can use for instance any $\varphi\in C^{\infty}(\R^{N})$,
supported in
${\cal{C}}=\{\xi\in\R^{N}/\frac{3}{4}\leq|\xi|\leq\frac{8}{3}\}$
such that:
$$\sum_{l\in\mathbb{Z}}\varphi(2^{-l}\xi)=1\,\,\,\,\mbox{if}\,\,\,\,\xi\ne 0.$$
Denoting $h={\cal{F}}^{-1}\varphi$, we then define the dyadic
blocks by:
$$\D_{l}u=\varphi(2^{-l}D)u=2^{lN}\int_{\R^{N}}h(2^{l}y)u(x-y)dy\,\,\,\,\mbox{and}\,\,\,S_{l}u=\sum_{k\leq
l-1}\D_{k}u\,.$$ Formally, one can write that:
$$u=\sum_{k\in\mathbb{Z}}\D_{k}u\,.$$
This decomposition is called homogeneous Littlewood-Paley
decomposition. Let us observe that the above formal equality does
not hold in ${\cal{S}}^{'}(\R^{N})$ for two reasons:
\begin{enumerate}
\item The right hand-side does not necessarily converge in
${\cal{S}}^{'}(\R^{N})$.
\item Even if it does, the equality is not
always true in ${\cal{S}}^{'}(\R^{N})$ (consider the case of the
polynomials).
\end{enumerate}
\subsection{Homogeneous Besov spaces and first properties}
\begin{definition}
For
$s\in\R,\,\,p\in[1,+\infty],\,\,q\in[1,+\infty],\,\,\mbox{and}\,\,u\in{\cal{S}}^{'}(\R^{N})$
we set:
$$\|u\|_{B^{s}_{p,q}}=(\sum_{l\in\mathbb{Z}}(2^{ls}\|\D_{l}u\|_{L^{p}})^{q})^{\frac{1}{q}}.$$
The Besov space $B^{s}_{p,q}$ is the set of temperate distribution $u$ such that $\|u\|_{B^{s}_{p,q}}<+\infty$.
\end{definition}
\begin{remarka}The above definition is a natural generalization of the
nonhomogeneous Sobolev and H$\ddot{\mbox{o}}$lder spaces: one can show
that $B^{s}_{\infty,\infty}$ is the nonhomogeneous
H$\ddot{\mbox{o}}$lder space $C^{s}$ and that $B^{s}_{2,2}$ is
the nonhomogeneous space $H^{s}$.
\end{remarka}
\begin{proposition}
\label{derivation,interpolation}
The following properties holds:
\begin{enumerate}
\item there exists a constant universal $C$
such that:\\
$C^{-1}\|u\|_{B^{s}_{p,r}}\leq\|\n u\|_{B^{s-1}_{p,r}}\leq
C\|u\|_{B^{s}_{p,r}}.$
\item If
$p_{1}<p_{2}$ and $r_{1}\leq r_{2}$ then $B^{s}_{p_{1},r_{1}}\hookrightarrow
B^{s-N(1/p_{1}-1/p_{2})}_{p_{2},r_{2}}$.
\item $B^{s^{'}}_{p,r_{1}}\hookrightarrow B^{s}_{p,r}$ if $s^{'}> s$ or if $s=s^{'}$ and $r_{1}\leq r$.
\end{enumerate}
\label{interpolation}
\end{proposition}
Let now recall a few product laws in Besov spaces coming directly from the paradifferential calculus of J-M. Bony
(see \cite{BJM}) and rewrite on a generalized form in \cite{AP} by H. Abidi and M. Paicu (in this article the results are written
in the case of homogeneous sapces but it can easily generalize for the nonhomogeneous Besov spaces).
\begin{proposition}
\label{produit1}
We have the following laws of product:
\begin{itemize}
\item For all $s\in\R$, $(p,r)\in[1,+\infty]^{2}$ we have:
\begin{equation}
\|uv\|_{B^{s}_{p,r}}\leq
C(\|u\|_{L^{\infty}}\|v\|_{B^{s}_{p,r}}+\|v\|_{L^{\infty}}\|u\|_{B^{s}_{p,r}})\,.
\label{2.2}
\end{equation}
\item Let $(p,p_{1},p_{2},r,\lambda_{1},\lambda_{2})\in[1,+\infty]^{2}$ such that:$\frac{1}{p}\leq\frac{1}{p_{1}}+\frac{1}{p_{2}}$,
$p_{1}\leq\lambda_{2}$, $p_{2}\leq\lambda_{1}$, $\frac{1}{p}\leq\frac{1}{p_{1}}+\frac{1}{\lambda_{1}}$ and
$\frac{1}{p}\leq\frac{1}{p_{2}}+\frac{1}{\lambda_{2}}$. We have then the following inequalities:\\
if $s_{1}+s_{2}+N\inf(0,1-\frac{1}{p_{1}}-\frac{1}{p_{2}})>0$, $s_{1}+\frac{N}{\lambda_{2}}<\frac{N}{p_{1}}$ and
$s_{2}+\frac{N}{\lambda_{1}}<\frac{N}{p_{2}}$ then:
\begin{equation}
\|uv\|_{B^{s_{1}+s_{2}-N(\frac{1}{p_{1}}+\frac{1}{p_{2}}-\frac{1}{p})}_{p,r}}\lesssim\|u\|_{B^{s_{1}}_{p_{1},r}}
\|v\|_{B^{s_{2}}_{p_{2},\infty}},
\label{2.3}
\end{equation}
when $s_{1}+\frac{N}{\lambda_{2}}=\frac{N}{p_{1}}$ (resp $s_{2}+\frac{N}{\lambda_{1}}=\frac{N}{p_{2}}$) we replace
$\|u\|_{B^{s_{1}}_{p_{1},r}}\|v\|_{B^{s_{2}}_{p_{2},\infty}}$ (resp $\|v\|_{B^{s_{2}}_{p_{2},\infty}}$) by
$\|u\|_{B^{s_{1}}_{p_{1},1}}\|v\|_{B^{s_{2}}_{p_{2},r}}$ (resp $\|v\|_{B^{s_{2}}_{p_{2},\infty}\cap L^{\infty}}$),
if $s_{1}+\frac{N}{\lambda_{2}}=\frac{N}{p_{1}}$ and $s_{2}+\frac{N}{\lambda_{1}}=\frac{N}{p_{2}}$ we take $r=1$.
\\
If $s_{1}+s_{2}=0$, $s_{1}\in(\frac{N}{\lambda_{1}}-\frac{N}{p_{2}},\frac{N}{p_{1}}-\frac{N}{\lambda_{2}}]$ and
$\frac{1}{p_{1}}+\frac{1}{p_{2}}\leq 1$ then:
\begin{equation}
\|uv\|_{B^{-N(\frac{1}{p_{1}}+\frac{1}{p_{2}}-\frac{1}{p})}_{p,\infty}}\lesssim\|u\|_{B^{s_{1}}_{p_{1},1}}
\|v\|_{B^{s_{2}}_{p_{2},\infty}}.
\label{2.4}
\end{equation}
If $|s|<\NN$ for $p\geq2$ and $-\frac{N}{p^{'}}<s<\NN$ else, we have:
\begin{equation}
\|uv\|_{B^{s}_{p,r}}\leq C\|u\|_{B^{s}_{p,r}}\|v\|_{B^{\NN}_{p,\infty}\cap L^{\infty}}.
\label{2.5}
\end{equation}
\end{itemize}
\end{proposition}
\begin{remarka}
In the sequel $p$ will be either $p_{1}$ or $p_{2}$ and in this case $\frac{1}{\lambda}=\frac{1}{p_{1}}-\frac{1}{p_{2}}$
if $p_{1}\leq p_{2}$, resp $\frac{1}{\lambda}=\frac{1}{p_{2}}-\frac{1}{p_{1}}$
if $p_{2}\leq p_{1}$.
\end{remarka}
\begin{corollaire}
\label{produit2}
Let $r\in [1,+\infty]$, $1\leq p\leq p_{1}\leq +\infty$ and $s$ such that:
\begin{itemize}
\item $s\in(-\frac{N}{p_{1}},\frac{N}{p_{1}})$ if $\frac{1}{p}+\frac{1}{p_{1}}\leq 1$,
\item $s\in(-\frac{N}{p_{1}}+N(\frac{1}{p}+\frac{1}{p_{1}}-1),\frac{N}{p_{1}})$ if $\frac{1}{p}+\frac{1}{p_{1}}> 1$,
\end{itemize}
then we have if $u\in B^{s}_{p,r}$ and $v\in B^{\frac{N}{p_{1}}}_{p_{1},\infty}\cap L^{\infty}$:
$$\|uv\|_{B^{s}_{p,r}}\leq C\|u\|_{B^{s}_{p,r}}\|v\|_{B^{\frac{N}{p_{1}}}_{p_{1},\infty}\cap L^{\infty}}.$$
\end{corollaire}
The study of non stationary PDE's requires space of type $L^{\rho}(0,T,X)$ for appropriate Banach spaces $X$. In our case, we
expect $X$ to be a Besov space, so that it is natural to localize the equation through Littlewood-Payley decomposition. But, in doing so, we obtain
bounds in spaces which are not type $L^{\rho}(0,T,X)$ (except if $r=p$).
We are now going to
define the spaces of Chemin-Lerner in which we will work, which are
a refinement of the spaces
$L_{T}^{\rho}(B^{s}_{p,r})$.
$\hspace{15cm}$
\begin{definition}
Let $\rho\in[1,+\infty]$, $T\in[1,+\infty]$ and $s_{1}\in\R$. We set:
$$\|u\|_{\widetilde{L}^{\rho}_{T}(B^{s_{1}}_{p,r})}=
\big(\sum_{l\in\mathbb{Z}}2^{lrs_{1}}\|\D_{l}u(t)\|_{L^{\rho}(L^{p})}^{r}\big)^{\frac{1}{r}}\,.$$
We then define the space $\widetilde{L}^{\rho}_{T}(B^{s_{1}}_{p,r})$ as the set of temperate distribution $u$ over
$(0,T)\times\R^{N}$ such that 
$\|u\|_{\widetilde{L}^{\rho}_{T}(B^{s_{1}}_{p,r})}<+\infty$.
\end{definition}
We set $\widetilde{C}_{T}(\widetilde{B}^{s_{1}}_{p,r})=\widetilde{L}^{\infty}_{T}(\widetilde{B}^{s_{1}}_{p,r})\cap
{\cal C}([0,T],B^{s_{1}}_{p,r})$.
Let us emphasize that, according to Minkowski inequality, we have:
$$\|u\|_{\widetilde{L}^{\rho}_{T}(B^{s_{1}}_{p,r})}\leq\|u\|_{L^{\rho}_{T}(B^{s_{1}}_{p,r})}\;\;\mbox{if}\;\;r\geq\rho
,\;\;\;\|u\|_{\widetilde{L}^{\rho}_{T}(B^{s_{1}}_{p,r})}\geq\|u\|_{L^{\rho}_{T}(B^{s_{1}}_{p,r})}\;\;\mbox{if}\;\;r\leq\rho
.$$
\begin{remarka}
It is easy to generalize proposition \ref{produit1},
to $\widetilde{L}^{\rho}_{T}(B^{s_{1}}_{p,r})$ spaces. The indices $s_{1}$, $p$, $r$
behave just as in the stationary case whereas the time exponent $\rho$ behaves according to H\"older inequality.
\end{remarka}
In the sequel we will need of composition lemma in $\widetilde{L}^{\rho}_{T}(B^{s}_{p,r})$ spaces.
\begin{lemme}
\label{composition}
Let $s>0$, $(p,r)\in[1,+\infty]$ and $u\in \widetilde{L}^{\rho}_{T}(B^{s}_{p,r})\cap L^{\infty}_{T}(L^{\infty})$.
\begin{enumerate}
 \item Let $F\in W_{loc}^{[s]+2,\infty}(\R^{N})$ such that $F(0)=0$. Then $F(u)\in \widetilde{L}^{\rho}_{T}(B^{s}_{p,r})$. More precisely there exists a function $C$ depending only on $s$, $p$, $r$, $N$ and $F$ such that:
$$\|F(u)\|_{\widetilde{L}^{\rho}_{T}(B^{s}_{p,r})}\leq C(\|u\|_{L^{\infty}_{T}(L^{\infty})})\|u\|_{\widetilde{L}^{\rho}_{T}(B^{s}_{p,r})}.$$
\item If $v,\,u\in\widetilde{L}^{\rho}_{T}(B^{s}_{p,r})\cap
L^{\infty}_{T}(L^{\infty})$ and $G\in
W^{[s]+3,\infty}_{loc}(\R^{N})$ then $G(u)-G(v)$ belongs to
$\widetilde{L}^{\rho}_{T}(B^{s}_{p})$ and there exists a constant C
depending only of $s, p ,N\;\mbox{and}\;G$ such that:
$$
\begin{aligned}
\|G(u)-G(v)\|_{\widetilde{L}^{\rho}_{T}(B^{s}_{p,r})}\leq&
\,\,C(\|u\|_{L^{\infty}_{T}(L^{\infty})},\|v\|_{L^{\infty}_{T}(L^{\infty})})
\big(\|v-u\|_{\widetilde{L}^{\rho}_{T}(B^{s}_{p,r})}
(1+\|u\|_{L^{\infty}_{T}(L^{\infty})}\\
&+\|v\|_{L^{\infty}_{T}(L^{\infty})})+\|v-u\|_{L^{\infty}_{T}(L^{\infty})}(\|u\|_{\widetilde{L}^{\rho}_{T}(B^{s}_{p,r})}
+\|v\|_{\widetilde{L}^{\rho}_{T}(B^{s}_{p,r})})\big).
\end{aligned}
$$
\end{enumerate}
\end{lemme}
Now we give some result on the behavior of the Besov spaces via some pseudodifferential operator (see \cite{BCD}).
\begin{definition}
Let $m\in\R$. A smooth function function $f:\R^{N}\rightarrow\R$ is said to be a ${\cal S}^{m}$ multiplier if for all muti-index $\alpha$, there exists a constant $C_{\alpha}$ such that:
$$\forall\xi\in\R^{N},\;\;|\p^{\alpha}f(\xi)|\leq C_{\alpha}(1+|\xi|)^{m-|\alpha|}.$$
\label{smoothf}
\end{definition}
\begin{proposition}
Let $m\in\R$ and $f$ be a ${\cal S}^{m}$ multiplier. Then for all $s\in\R$ and $1\leq p,r\leq+\infty$ the operator $f(D)$ is continuous from $B^{s}_{p,r}$ to $B^{s-m}_{p,r}$.
\label{singuliere}
\end{proposition}
Let us now give some estimates for the heat equation:
\begin{proposition}
\label{5chaleur} Let $s\in\R$, $(p,r)\in[1,+\infty]^{2}$ and
$1\leq\rho_{2}\leq\rho_{1}\leq+\infty$. Assume that $u_{0}\in B^{s}_{p,r}$ and $f\in\widetilde{L}^{\rho_{2}}_{T}
(B^{s-2+2/\rho_{2}}_{p,r})$.
Let u be a solution of:
$$
\begin{cases}
\begin{aligned}
&\p_{t}u-\mu\D u=f\\
&u_{t=0}=u_{0}\,.
\end{aligned}
\end{cases}
$$
Then there exists $C>0$ depending only on $N,\mu,\rho_{1}$ and
$\rho_{2}$ such that:
$$\|u\|_{\widetilde{L}^{\rho_{1}}_{T}(\widetilde{B}^{s+2/\rho_{1}}_{p,r})}\leq C\big(
 \|u_{0}\|_{B^{s}_{p,r}}+\mu^{\frac{1}{\rho_{2}}-1}\|f\|_{\widetilde{L}^{\rho_{2}}_{T}
 (B^{s-2+2/\rho_{2}}_{p,r})}\big)\,.$$
 If in addition $r$ is finite then $u$ belongs to $C([0,T],B^{s}_{p,r})$.
\end{proposition}
\section{Proof of the theorem \ref{theo1}}
\label{section3}
The existence part of the theorem is proved by an iterative method. We define a sequence $(q^{n},m^{n})$ such that:
$$
\begin{cases}
\begin{aligned}
&\p_{t}q^{0}-\mu\D q^{0}+{\ rm div}m^{0}=0,\\
&\p_{t}m^{0}-\mu\D m^{0}+r m^{0}=0,\\
&(q^{0},m^{0})=(q_{0},m_{0}).
\end{aligned}
\end{cases}
$$
Assuming that $(q^{n},m^{n})$ is in $E_{T}$ with:
$$E_{T}=\big(\widetilde{C}_{T}(B^{\N}_{2,1})\cap L^{1}_{T}(B^{\N+2}_{2,1})\big)\times\big(
\widetilde{C}_{T}(B^{\N-1}_{2,1})\cap L^{1}_{T}(B^{\N+1}_{2,1}\cap B^{\N-1}_{2,1})\big)^{N},$$
we define then  $q^{n+1}=q^{0}+\bar{q}^{n+1}$, $m^{n+1}=m^{0}+\bar{m}^{n+1}$  such that 
$(\bar{q}_{n+1},\bar{m}_{n+1})$ be the solution
of the following system:
$$
\begin{cases}
\begin{aligned}
&\p_{t}\bar{q}^{n+1}-\mu\D\bar{q}^{n+1}+{\rm div}\bar{m}^{n+1}=0,\\
& \p_{t}\bar{m}^{n+1}-\mu\D\bar{m}^{n+1}+r \bar{m}^{n+1} =G^{n},\\
&(\bar{q}^{n+1},\bar{m}^{n+1})_{/t=0}=(0,0),
\end{aligned}
\end{cases}
$$
with:
$$
\begin{aligned}
G^{n}=&-{\rm div}(\frac{m^{n}}{h^{n}}\otimes m^{n})
\end{aligned}
$$
We also set: $h^{n}=q^{n}+1$.
\subsubsection*{1) First step, uniform bounds:}
 Let $\e$ be a small
parameter and by proposition \ref{5chaleur}, we have for any $T>0$:
\begin{equation}
\begin{aligned}
&\|q^{0}\|_{L^{\infty}_{T}(B^{\N}_{2,1})\cap L^{1}_{T}(B^{\N+2}_{2,1})}\leq C\|q_{0}\|_{B^{\N}_{2,1}},\\
&\|m^{0}\|_{L^{\infty}_{T}(B^{\N-1}_{2,1})\cap L^{1}_{T}(B^{\N-1}_{2,1}\cap B^{\N+1}_{2,1})}\leq C\|m_{0}\|_{B^{\N-1}_{2,1}}.
\end{aligned}
\label{initial}
\end{equation}
We are going to show by induction that for $\e>0$ small enough:
$$\|(\bar{q}^{n},\bar{m}^{n})\|_{F_{T}}\leq\e.\leqno{({\cal{P}}_{n})}$$
As $(\bar{q}_{0},\bar{m}_{0})=(0,0)$ the result
is true for $n=0$. We suppose now $({\cal{P}}_{n})$ true and we are
going to show $({\cal{P}}_{n+1})$.
\\
To begin with we are going to show that $1+q^{n}$ is  positive. Indeed we have: $h^{0}=h^{0}_{1}+h^{0}_{2}$ such that:
$$
\begin{aligned}
&\p_{t}h^{0}_{1}-\mu \D h^{0}_{1}=0,\\
&(h^{0}_{1})_{/t=0}=h_{0}.
\end{aligned}
$$
and:
$$
\begin{aligned}
&\p_{t}h^{0}_{2}-\mu \D h^{0}_{2}=-{\rm div}m^{0},\\
&(h^{0}_{1})_{/t=0}=h_{0}.
\end{aligned}
$$
By proposition (\ref{5chaleur}) and (\ref{initial}) we have for any $T>0$:
\begin{equation}
\|h^{0}_{2}\|_{L^{\infty}_{T}(B^{\N}_{2,1})}\leq C\|m_{0}\|_{B^{\N-1}_{2,1}}.
\label{ss5}
\end{equation}
By maximum principle, we have for any $t>0$:
$$h^{0}_{1}(t,x)\geq \min_{x\in\R^{N}}h_{0}(x)\geq c>0.$$
We deduce that for $\eta=\|m_{0}\|_{B^{\N-1}_{2,1}}$ (at least inferior to $\frac{c}{4C}$ with the $C$ of (\ref{ss5})) small enough and any $t>0$:
$$h^{0}(t,x)/geq \frac{3c}{4}>0,$$
and
$$q^{0}(t,x)\geq \frac{3c}{4}-1.$$
and by definition of $q^{n}$ and the assumption ${\cap P}_{n}$ that:
$$q^{n}(t,x)\geq  \frac{3c}{4}-1-\e.$$
In particular for $\e$ small enough at least $\e\leq\frac{c}{4}$, we deduce that:
\begin{equation}
h^{n}=1+q^{n}\geq \frac{c}{2}>0.
\label{vide}
\end{equation}
In order to bound $(\bar{q}^{n},\bar{m}^{n})$ in $E_{T}$, we shall use proposition \ref{5chaleur} and in particular estimating 
$G^{n}$ 
in $L^{1}_{T}(B^{\N-1}_{2,1})$.
By using proposition \ref{produit1}, (\ref{vide}) and lemma\ref{composition}, we obtain:
\begin{equation}
\begin{aligned}
&\|{\rm div}(\frac{m^{n}}{h^{n}}\otimes m^{n})\|_{L^{1}_{T}(B^{\N-1}_{2,1})}\leq \|\frac{m^{n}}{h^{n}}\otimes m^{n}\|_{L^{1}_{T}(B^{\N-1}_{2,1})},\\
&\leq C\|m^{n}\|_{L^{2}_{T}(B^{\N}_{2,1})}^{2}\big(\|\frac{1}{1+q^{n}}-1\|_{L^{\infty}_{T}(B^{\N}_{2,1})}+1\big),\\
&\leq C(\|m^{0}\|_{L^{2}_{T}(B^{\N}_{2,1})}^{2}+\|\bar{m}^{n}\|_{L^{2}_{T}(B^{\N}_{2,1})})^{2}(1+C(\|\frac{1}{h^{n}}\|_{L_{T}^{\infty}}))\big(\|q^{n}\|_{L_{T}^{\infty}(B^{\N}_{2,1})}+1\big),\\
&\leq C(\|m^{0}\|_{L^{2}_{T}(B^{\N}_{2,1})}^{2}+\|\bar{m}^{n}\|_{L^{2}_{T}(B^{\N}_{2,1})})^{2}\big(\|\bar{q}^{n}\|_{L_{T}^{\infty}(B^{\N}_{2,1})}+\|q^{0}\|_{L_{T}^{\infty}(B^{\N}_{2,1})}
+1\big),
\end{aligned}
\label{a2}
\end{equation}
Therefore by using 
(\ref{a2}), the proposition \ref{5chaleur} and $({\cal P}_{n})$ we obtain for any $T>0$:
\begin{equation}
\begin{aligned}
\|(\bar{q}^{n+1},\bar{m}^{n+1})\|_{F_{T}}&\leq C (\|m^{0}\|_{L^{2}_{T}(B^{\N}_{2,1})}^{2}+\e)^{2}\big(\e+\|q^{0}\|_{L_{T}^{\infty}(B^{\N}_{2,1})}
+1\big),\\
&\leq C (\eta+\e)^{2}\big(2+\|q^{0}\|_{L_{T}^{\infty}(B^{\N}_{2,1})}\big)
\end{aligned}
\end{equation}
By choosing $\eta=\e$ and $\e\leq \frac{1}{2C(2+\|q^{0}\|_{L_{T}^{\infty}(B^{\N}_{2,1})})}$, this implies $({\cal P})_{n+1}$.  We
have shown by induction that $(q^{n},m^{n})$ is uniformly bounded
in $F_{T}$ for any $T>0$.
\subsubsection*{Second Step: Convergence of the
sequence}
 We shall prove
that $(q^{n},m^{n})$ is a Cauchy sequence in the Banach
space $F_{T}$, hence converges to some
$(q,m)\in F_{T}$.\\
Let:
$$\delta q^{n}=q^{n+1}-q^{n}\;\;\mbox{and}\;\;\delta m^{n}=m^{n+1}-m^{n}.$$
The system verified by $(\de q^{n},\de m^{n})$ reads:
$$
\begin{cases}
\begin{aligned}
&\p_{t}\delta q^{n}-\mu\D \delta q^{n}+{\rm div}\delta m^{n}=0,\\
&\p_{t}\delta m^{n}-\mu\D\delta m^{n}+r \delta m^{n} =G^{n}-G{n-1},\\
&\delta q^{n}(0)=0\;,\;\delta u^{n}(0)=0.
\end{aligned}
\end{cases}
$$
Applying propositions \ref{5chaleur} and using $({\cal{P}}_{n})$, we get for any $T>0$:
$$
\begin{aligned}
&\|(\de q^{n},\de m^{n})\|_{F_{T}}\leq\;
C\|G^{n}-G_{n-1}\|_{L^{1}_{T}(B^{N/2-1})},\\
&\leq C\big(\|\frac{\delta m^{n}}{h^{n}}\otimes m^{n}\|_{L^{1}_{T}(B^{N/2}_{2,1})}+
\|\frac{\delta m^{n}}{h^{n}}\otimes m^{n-1}\|_{L^{1}_{T}(B^{N/2}_{2,1})}+\|m^{n}\otimes m^{n-1}(\frac{1}{h^{n}}-\frac{1}{h^{n-1}})\|_{L^{1}_{T}(B^{N/2}_{2,1})}\big).
\end{aligned}
$$
By using proposition \ref{produit1} and lemma \ref{composition}, we get:
$$\|(\de q^{n},\de m^{n})\|_{F_{T}}\leq
C\e\|(\de q^{n-1},\de m^{n-1})\|_{F_{T}}.$$ So by taking $\e$ enough small we have proved that
$(q^{n},m^{n})$ is a Cauchy sequence in $F_{T}$ which is a Banach space. It implies that $(q^{n},m^{n})$ converge to a limit $(q,m)$  in $F_{T}$.  It is easy to verify that $(q,m)$ is a
solution of the system (\ref{3system4}).
\subsubsection*{3)Uniqueness of the solution:}
The proof is similar to the proof of contraction, indeed we need the same type
 of estimates. Let us consider two solutions in
$E_{T}$: $(q_{1},m_{1})$ and
$(q_{2},m_{2})$ of the system (\ref{3system4}) with the same
initial data. With no loss of generality, one can assume that
$(q_{1},m_{1})$ is the solution found in the previous
section. We thus have:
$$q_{1}(t,x)\geq -\frac{1}{2}.
\leqno{\cal{(H)}}$$ 
We note:
$$\delta q=q_{2}-q_{1},\;\delta m=m_{2}-m_{1},$$ 
which verifies the system:
$$
\begin{cases}
\begin{aligned}
&\p_{t}\delta q-\mu\D\delta q+{\rm div}\delta m=0,\\
&\p_{t}\delta m-\mu\D\delta m+r \delta m=-{\rm div}(\frac{m_{1}}{h_{1}}\otimes m_{1})+
{\rm div}(\frac{m_{1}}{h_{1}}\otimes m_{1})
\end{aligned}
\end{cases}
$$
By using proposition \ref{produit1}, \ref{5chaleur} and lemma \ref{composition} on $[0,T_{1}]$
with $0<T$ we have:
$$\|(\delta q,\delta m)\|_{E_{T}}\leq A(T)\|(\delta q,\delta m)\|_{E_{T}},$$
such that for $T$ small enough $A(T)\leq\frac{1}{2}$. We thus obtain: $\de q=0$, $\de m=0$ on
$[0,T]$. And we repeat the argument in order to prove that:
$\delta q=0$, $\de m=0$ on $\R^{+}$. This conclude the proof of theorem \ref{theo1}.
\hfill {$\Box$}

\end{document}